\documentclass{article}

\usepackage{url}

\usepackage{amsmath, amssymb, amstext, amsopn, amsxtra, graphicx, color, calligra, hyperref, url, ulem, polynom, permute, ifthen, float, epstopdf, linsys, marvosym, rotating}

\begin{document}

\newcommand{\recip}[1]{\ensuremath{\frac1{#1}}}
\newcommand{\half}{\recip2}
\newcommand{\third}{\recip{3}}
\newcommand{\twothirds}{\ensuremath{\frac23}}
\newcommand{\fourth}{\recip4}
\newcommand{\fifth}{\recip5}
\newcommand{\sixth}{\recip6}
\newcommand{\point}[1]{\mbox{\(\left(#1\right)\)}}
\newcommand{\qar}{\begin{eqnarray*}}
\newcommand{\raq}{\end{eqnarray*}}
\newcommand{\inv}{^{-1}}
\newcommand{\vone}{{\bf 1}}
\newcommand{\va}{{\bf a}}
\newcommand{\vb}{{\bf b}}
\newcommand{\vr}{{\bf r}}
\newcommand{\vs}{{\bf s}}
\newcommand{\vu}{{\bf u}}
\newcommand{\comment}[1]{\hskip 5mm \parbox{3in}{#1}}
\newcommand{\casI}{\renewcommand{\labelenumi}{\bf Case \Roman{enumi}.}\num}
\newcommand{\sacI}{\mun \renewcommand{\labelenumi}{\arabic{enumi}.}}
\newcommand{\num}{\begin{enumerate}}
\newcommand{\mun}{\end{enumerate}}
\newcommand{\rowr}{\ensuremath{\left(\begin{array}{rrrr}r_1 & r_2 & \cdots & r_k \end{array}\right)}}
\newcommand{\rows}{\ensuremath{\left(\begin{array}{rrrr}s_1 & s_2 & \cdots & s_k \end{array}\right)}}
\newcommand{\qed}{\mbox{\(\square\,\,\)}}
\newcommand{\longcomment}[1]{}

\newtheorem{thm}{Theorem}
\newtheorem{cor}[thm]{Corollary}
\newtheorem{rem}[thm]{Remark}

\renewcommand{\em}{\textit}

\title{Markov Chains for Collaboration}

\author{Robert Mena\\
rmena@csulb.edu\\
Will Murray\\
wmurray@csulb.edu\\
California State University, Long Beach\\
Long Beach, CA 98040-1001}
\date{January 13, 2014}

\maketitle

\section{Introduction:  Who wants to be a collaborator?}

The math department at New Alarkania State University is comprised of Alan the analyst, Lorraine the logician, Stacy the statistician, and Tom the topologist.  Each one is desperate for collaborators, so they start a Friday poker series.  Each one is equally skilled, and they agree that the loser of each week's game (the first to run out of money) will renounce his or her former field and join the research team of the biggest winner.  

In the first week, Stacy wins and Tom loses, so Tom gives up topology and joins Stacy to study statistics.  The following week, Lorraine wins and Stacy loses, so Stacy becomes a logician.  Next, Stacy wins and Lorraine loses, so no one has to switch.  You have no doubt already guessed that eventually (with probability one), all of them will be working in the same field.  (After the first week, for example, one field has already disappeared permanently, since as soon as Tom loses there are no more topologists.)  

This is an example of a {\em {Markov chain}}, in which a system can be in a number of possible states, and at each time step there is a certain probability of moving to each of the other states (or of remaining in the same state).  Kemeny and Snell (\cite{kemeny}) give an excellent background on Markov chains.

We will break our chain up into {\em {stages}}, numbered in reverse order according to how many fields are remaining.  Thus, we start in Stage 4, meaning there are four fields left, but after one week we are certain to be in Stage 3.  We will study three questions here:

\begin{enumerate}
\item  \label{question-stagetime} {\em {How long do we expect to stay in each stage?}}  The expected time in Stage 4 (or Stage \(n\) in the general case of \(n\) starters) is exactly one week, but after that it gets more complicated.
\item  \label{question-landvector} {\em {When we first arrive at Stage \(t-1\) from Stage \(t\), what is the most likely configuration of the fields?}}  More precisely, what are the probablilities of arriving at different configurations of the players into \(t-1\) teams?  For example, with \(n=4\) starters, when we go down from three fields to two, are we more likely to have two teams of two players each, or a team of three and a lone wolf?
\item  \label{question-totaltime} {\em {How long does the game last?}}  In other words, what is the expected time until we reach the absorbing state in which everyone is on the same team?  Of course, the answer here is just the sum of the answers from Question~\ref{question-stagetime}.
\end{enumerate}

We invite you to play with the small cases of \(n = 3, 4\), or 5 starters, which are not too hard to work out from first principles.  You will find that the answers to Question~\ref{question-totaltime} are 4, 9, and 16 weeks respectively.  It might not be obvious that this stunning pattern should continue to hold, but we will prove that with \(n\) starters, the expected time is indeed \((n-1)^2\) weeks.  (Unfortunately, there appears to be no correspondingly congenial answer for the variance.)

The general answers to Questions~\ref{question-stagetime} and~\ref{question-landvector} are not so obvious from analyzing small cases.  For example, with \(n=5\) starters, the total expected time of 16 weeks breaks down into stages of \(e_{54} = 1, e_{43} = \frac{5}{3}, e_{32} = \frac{10}{3}\), and  \(e_{21} = 10\) weeks.  We will see that these come from binomial coefficients and that the answer to Question~\ref{question-landvector} comes from multinomial coefficients.


We organize the paper as follows:  In the second section, we will warm up by solving the case \(n=4\) from scratch, using no sophisticated machinery.  Besides resolving the question for New Alarkania State, this will give us an informal preview of some of the notation and theorems  coming later.  Next, we introduce more formal notation and illustrate it with a larger example, \(n=6\).  We then study the vectors of probabilities and discover multinomial coefficients as the answer to Question~\ref{question-landvector}.  With the probability vectors in hand, it is relatively quick to study the expected times and answer Questions~\ref{question-stagetime} and~\ref{question-totaltime}.  In the final section, we present a symmetric approach that answers Question~\ref{question-totaltime} directly without reference to the answers to Questions~\ref{question-stagetime} and~\ref{question-landvector}.

\section{\(n=4\):  How long must New Alarkania wait?}

In this section we will work out the case of four players from scratch using only basic probability; however, some of the notation and theory for later will become evident as we go along.  As mentioned above, we organize the possible configurations into {\em {stages}} according to the number of teams left; thus we proceed in reverse order from Stage 4 (four individuals, \([1111]\)) down to Stage 1 (a single team of four, \([4]\)).  

Starting at Stage 4 (\([1111]\)), note that in the first week, one player must lose and join the winner's team.  Therefore, the expected time to Stage 3 is exactly \(e_{43} = 1\) week.  The configuration at Stage 3 is necessarily \([211]\), one team of two players and two individuals.  

Now, from \([211]\), the loser can be one of the players on the team of two, in which case the new configuration is still \([211]\).  (If the winner is the other player on the team, then there is no change at all; if the winner is one of the two individuals, then the loser joins that individual, making a new team of two and leaving the loser's former teammate as an individual.)  If the loser is one of the two individuals, however, we will go down to Stage 2.  The new configuration depends on who the winner is, but we note first that since there is a \half\ chance of the loser being one of the two individuals, the expected waiting time is exactly \(e_{32} = 2\) weeks.  

When we do first get down to Stage 2, what configuration will we land in?  We know that the loser in the previous week was one of the two individuals.  There is a \twothirds\ chance that the winner was a member of the team of two, in which case we land in \([31]\).  There is a \third\ chance that the winner was the other individual, landing us in \([22]\).  We thus have an answer for Question~\ref{question-landvector} at Stage 2:  We say \(L_2 := \left(\begin{array}{rr}\twothirds & \third \end{array}\right)\) is the {\em {landing vector}} at Stage 2, representing the probabilities that when we first arrive in Stage 2, we land in \([31]\) or \([22]\) respectively.  (We had landing vectors at the previous stages as well, but because there was only one configuration in each stage, they were simply the trivial vectors \(L_4 := (1), L_3 := (1)\).)

Finally, we calculate the expected time \(e_{21}\) to go from Stage 2 to Stage 1.  Here are the possible outcomes from configuration \([31]\):

\bigskip

\begin{tabular}{cll}
Probability & Outcome & Explanation\\ \hline
\half & Stay at \([31]\). & Winner and loser are both from the team of three.\\
\fourth & Move to \([22]\). & Winner is the individual.\\
\fourth & Move to \([4]\). & Loser is the individual.  
\end{tabular}

\bigskip

And here are the possiblities from \([22]\):

\bigskip

\begin{tabular}{cll}
Probability & Outcome & Explanation\\ \hline
\twothirds & Move to \([31]\). & Winner and loser are from different teams.\\
\third & Stay at \([22]\). & Winner and loser are on the same team.\\
0 & Move to \([4]\). & Not possible in one week.
\end{tabular}

\bigskip

We collect these probabilities in a matrix, denoted \(A_2\), for later:

\bigskip

\[
\begin{array}{crccccrcc} & [31] & [22] & [4] & & & [31] & [22] & [4] \\
\begin{array}{c} [31] \\ \left[22\right] \end{array} & 
\left(\begin{array}{c} \half \\ \twothirds \end{array} \right.& 
\left.\begin{array}{c} \fourth \\ \third \end{array} \right.& 
\left|\begin{array}{r} \fourth \\ 0 \end{array} \right) & 
= & 
\left. \begin{array}{c} [31] \\ \left[22\right] \end{array} \right( & 
\multicolumn{2}{c}{A_2}& 
\left|\begin{array}{r} \fourth \\ 0 \end{array} \right)
\end{array}
\]

To find the expected time \(e_{21}\) to go from Stage 2 to Stage 1, let \(x_1\) be the expected time to go from \([31]\) to \([4]\) and let \(x_2\) be the expected time to go from \([22]\) to \([4]\) (necessarily via \([31]\)).  If we start at \([31]\) and let one week go by, there is a \half\ chance that we will stay at \([31]\), giving us a new expected time of \(x_1\) plus the one week that just elapsed.  There is a \fourth\ chance that we move to \([22]\), giving us a new expected time of \(x_2\) plus one.  Finally, there is a \fourth\ chance that we move directly to \([4]\), making the time exactly one week.  We summarize this as an equation:
\[
x_1 = \fourth \point{x_1+1} + \half \point{x_2+1} + \fourth \point{1} = \fourth x_1 + \half x_2 + 1
\]
Starting at \([22]\) and letting one week elapse gives us a similar equation:
\[
x_2 = \twothirds \point{x_1+1} + \third \point{x_2+1} + 0 \point{1} = \twothirds x_1 + \third x_2 + 1
\]

Combining these equations gives us a matrix equation that is easy to solve:
\qar
 \left(\begin{array}{c} x_1 \\ x_2 \end{array}\right) & = & A_2 \left(\begin{array}{c} x_1 \\ x_2 \end{array}\right) + \left(\begin{array}{c} 1 \\ 1 \end{array}\right)\\
\point{I - A_2}\left(\begin{array}{c} x_1 \\ x_2 \end{array}\right) & = & \left(\begin{array}{c} 1 \\ 1 \end{array}\right)\\
\left(\begin{array}{c} x_1 \\ x_2 \end{array}\right) & = & \point{I - A_2}\inv \left(\begin{array}{c} 1 \\ 1 \end{array}\right)\\
& = & \left(\begin{array}{c} \frac{11}2 \\ 7 \end{array}\right)
\raq
Recalling the landing vector of probabilities that we arrive at Stage 2 either in \([31]\) or \([22]\), the expected time to go to Stage 1 is then
\[
e_{21} = \left(\begin{array}{rr}\twothirds & \third \end{array}\right)\left(\begin{array}{c} \frac{11}2 \\ 7 \end{array}\right) = 6 \mbox{ weeks.}
\]
Finally, the total expected time to go from Stage 4 down to Stage 1 is the sum of the expected times at each stage, \(e_{43} + e_{32} + e_{21} = 1 + 2 + 6 = 9\) weeks, or \((n-1)^2\) for \(n = 4\).  

Besides answering Questions~\ref{question-stagetime} -~\ref{question-totaltime} for New Alarkania, this small example already showcases several features that will be reflected in larger cases later:  
\begin{itemize}
\item  We depended heavily on {\em {linearity of expectation}} to break the total expected time into a sum of expected times \(e_{t,t-1}\) to go from each Stage \(t\) to Stage \(t-1\).  

\item  Stage 2 (and for larger cases, almost all stages) consisted of multiple possible configurations, \([31]\) and \([22]\).  We described our arrival at Stage 2 in terms of a {\em {landing vector}} of probabilities \(L_2 := \left(\begin{array}{rr}\twothirds & \third \end{array}\right)\) that we would first land in each configuration.  These landing vectors are the answer to Question~\ref{question-landvector}, but this one small example is not enough to see the general pattern.

\item  We can compute the expected time to go from Stage \(t\) to Stage \(t-1\) as
\[
e_{t,t-1} = L_t(I-A_t)\inv \vone,\label{equation-stagetime}
\]
where \(L_t\) is the landing vector of probabilities for the configurations in Stage \(t\), \(A_t\) is the matrix of internal transition probabilities between the various configurations in Stage \(t\), and \(\vone\) is a column vector of ones of the appropriate length.

\item  In this small example, the expected times were all integers, \(e_{43} = 1, e_{32} = 2\), and \(e_{21} =  6\).  That won't generalize, but they will follow a most interesting pattern.  (We invite you to guess it now, with the reminder that
the times for the case \(n=5\) are \(e_{54} = 1, e_{43} = \frac{5}{3}, e_{32} = \frac{10}{3}\), and  \(e_{21} = 10\), giving a total time of \(1 + \frac{5}{3} + \frac{10}{3} + 10 = 16 = (n-1)^2\) weeks.)

\end{itemize}

Keeping the lessons from \(n=4\) in mind, we now move on to address the general problem.  

\section{Notation and examples}

Fix a value of \(n\).  We will consider the various partitions of \(n\) to be the states of the system.  We will use both {\em {partition notation}}, where we list the parts as \(n_1+n_2+\cdots +n_k\), which we will abbreviate as \(n_1n_2\cdots n_k\), 
 and {\em {vector notation}}, where we list the number of parts of each size as \((r_1r_2 \cdots r_k)\), so \(\sum i r_i = n\).  (When using vector notation, we will always assume that the last entry is nonzero.)

Let \(S(n,t)\) be the set of partitions of \(n\) in \(t\) parts, i.e. the set of all possible configurations at Stage \(t\).  Then the set of all partitions of \(n\) is \(\cup_{t=1}^n S(n,t)\).  We list the sets \(S(n,t)\) in reverse order from \(t=n\) to \(t=1\), and we assume that each \(S(n,t)\) is given a consistent internal ordering.

For example, let \(n=6\).  Then the states in partition notation are
\[
\{[111111],[21111],[2211,3111], [222,321,411],[33,42,51],[6]\},
\]
and, respectively, in vector notation are
\[
\{[(6)],[(41)],[(22),(301)],[(03),(111),(2001)],[(002),(0101),(10001)],[(000001)]\}.
\]

Let \(P\) be the probability transition matrix between the various possible states.  Then \(P\) is block upper bidiagonal, where each diagonal block is \(A_t\), the probability transition matrix from states in Stage \(t\) to each other, and each superdiagonal block is \(A_{t,t-1}\), the probability transition matrix from states in Stage \(t\) to states in Stage \(t-1\).  

For \(n=6\), using the ordering above, we have the following matrix:
\qar
P & = & \left(\begin{array}{l|l|l|l|l|l}
A_6 & A_{65} & & & &\\  \hline
& A_5 & A_{54} & & &\\  \hline
& & A_4 & A_{43} & &\\  \hline
& & & A_3 & A_{32} &\\  \hline
& & & & A_2 & A_{21}\\  \hline
& & & & & A_1
\end{array}\right)\\
& = &
\recip{30}
\left(\begin{array}{r|r|rr|rrr|rrr|r}
0 & 30 & & & & & & & & &\\  \hline
& 10 & 12 & 8 & & & & & & &\\  \hline
& & 12 & 8 & 2 & 8 & 0 & & & &\\
& & 9 & 6 & 0 & 6 & 9 & & & &\\  \hline
& & & & 6 & 24 & 0 & 0 & 0 & 0 &\\
& & & & 3 & 16 & 6 & 2 & 3 & 0 &\\
& & & & 0 & 8 & 12 & 0 & 2 & 8 &\\  \hline
& & & & & & & 12 & 18 & 0 & 0\\
& & & & & & & 8 & 14 & 8 & 0\\
& & & & & & & 0 & 5 & 20 & 5\\  \hline
& & & & & & & & & & 30
\end{array}\right)
\raq
For example, the middle rows of \(A_3\) and \(A_{32}\) are obtained by noting that of the 30 possible choices for winner and loser from the partition \(321\) (in vector notation, (111)), 3 lead to the partition \(222\), 16 to \(321\), 6 to \(411\), 2 to \(33\), 3 to \(42\), and none to \(51\)  (in vector notation, (03),(111),(2001),(002),(0101),(10001), respectively).  We invite you to check the other values.  

\section{Probability vectors and multinomial coefficients}

We define the {\em {landing vectors}} \(L_t\) recursively as follows.  First, we set \mbox{\(L_n := (1)\)} since we must start in Stage \(n\) in state \((n)\).  Now, for \(n \geq t \geq 2\), assume that we start in one of the states in Stage \(t\) with probabilities given by the entries of \(L_t\).  We then define \(L_{t-1}\) to be the row vector whose \(j\)-th entry is the probability that our first arrival in Stage \(t-1\) from Stage \(t\) is in the \(j\)-th state in Stage \(t-1\).

Thus, in the example above with \(n=6\), we have \(L_6 = (1), L_5 = (1)\), and \(L_4 = \left(\begin{array}{rr}\frac{3}{5} & \frac{2}{5}\end{array}\right)\), because when we move from Stage 5, necessarily starting at  \((41)\) (in vector notation), to Stage 4, 
we have a \(\frac{3}{5}\) chance of arriving in state (22) and a \(\frac{2}{5}\) chance of arriving in state (301).  

To calculate the \(L_t\)'s, we define \(P_{t,t-1}\) to be a matrix in which each row corresponds to a state in Stage \(t\) and each column to a state in Stage \(t-1\).  Entry \point{i,j} in \(P_{t,t-1}\) is defined to be the probability that, given that we start in state \(i\) in Stage \(t\), our first arrival in Stage \(t-1\) is in state \(j\).  By a similar derivation to the one we used in the example with \(n=4\) above, we have
\[
P_{t,t-1} = (I - A_t)\inv A_{t,t-1},
\]
where \(I\) is the identity matrix of appropriate size.  (This is also a standard result in the theory of Markov chains; see Theorem 3.3.7 in \cite{kemeny}.)

We can now compute the \(L_t\)'s recursively:
\[
L_{t-1} = L_tP_{t,t-1}= L_t (I - A_t)\inv A_{t,t-1}.  
\]

For example, with \(n=6\), we have 
\qar
L_5 & = & L_6 (I-A_6)\inv A_{65} = (1)(1)(1) = (1)\\
L_4 & = & L_5 (I-A_5)\inv A_{54} = (1)\left(\begin{array}{r}\frac{3}{2}\end{array}\right)\left(\begin{array}{rr}\frac{2}{5} & \frac{4}{15}\end{array}\right) = \left(\begin{array}{rr}\frac{3}{5} & \frac{2}{5}\end{array}\right)\\
L_3 & = & L_4 (I-A_4)\inv A_{43}\\
& = & \left(\begin{array}{rr}\frac{3}{5} & \frac{2}{5}\end{array}\right)\left(\begin{array}{rr}2 & \frac{2}{3}\\ \frac{3}{4} & \frac{3}{2}\end{array}\right)\left(\begin{array}{ccc}\frac{1}{15} & \frac{4}{15} & 0\\ 0 & \frac{1}{5} & \frac{3}{10}\end{array}\right)\\
& = & \left(\begin{array}{rrr}\frac{1}{10} & \frac{3}{5} & \frac{3}{10}\end{array}\right)
\raq
and so on.

We define \(VS(n,t)\) to be the vector space whose basis is the set of partitions \(S(n,t)\) in Stage \(t\).  Using vector notation for partitions, for
\[
\vr = \rowr \in S(n,t),
\]
we define the {\em {multinomial coefficient}}
\[
m_\vr := {t \choose r_1,r_2,\dots,r_k} = \frac{t!}{r_1!r_2!\cdots r_k!}.
\]
(Undergraduates will recall multinomial coefficients from combinatorial exercises about rearranging the letters of words like MISSISSIPPI; see Section 5.4 in \cite{brualdi} for details.)

Finally we define the vector \(\vu_t \in VS(n,t)\) by
\[
\vu_t := \sum_{\vr \in S(n,t)} m_\vr \vr
\]
and consider it as a row vector whose entries are the \(m_\vr\)'s.  

We can add the entries of a vector by multiplying it by \vone, the column vector of appropriate size whose entries are all ones.  

\begin{rem} \label{remark-sumcoeffs}The sum of the coefficients of \(\vu_t\) is 
\[
\vu_t \vone = \sum_{\vr \in S(n,t)} m_\vr = {n-1 \choose t-1}.
\]
\end{rem}

\textit{Proof}.  One way to list the partitions of \(n\) into \(t\) parts is to make a line of \(n\) pebbles and then insert \(t-1\) dividers into the \(n-1\) spaces between the pebbles; there are \({n-1 \choose t-1}\) ways to do this.  However, most partitions will be counted multiple times in this list since the parts can appear in any order.  In fact, the partition \(\vr = \rowr \in S(n,t)\) will appear exactly
\[
m_\vr = {t \choose r_1,r_2,\dots,r_k} = \frac{t!}{r_1!r_2!\cdots r_k!}
\]
times, giving the desired result.  \qed
\bigskip

For example, with \(n=6\), we have the following:
\[
\begin{array}{llllll}
\vu_6 & = & 1(6) & = & \left(\begin{array}{r}1\end{array}\right); & \vu_6\vone = 1\\
\vu_5 & = & 5(41)& = & \left(\begin{array}{r}5\end{array}\right); & \vu_5\vone = 5\\
\vu_4 & = & 6(22) + 4(301)& = & \left(\begin{array}{rr}6 & 4\end{array}\right); & \vu_4\vone = 10\\
\vu_3 & = & 1(03) + 6(111) + 3(2001)& = & \left(\begin{array}{rrr}1 & 6 & 3\end{array}\right); & \vu_3\vone = 10\\
\vu_2 & = & 1(002) + 2(0101) + 2(10001)& = & \left(\begin{array}{rrr}1 & 2 & 2\end{array}\right); & \vu_2\vone = 5\\
\vu_1 & = & 1(000001)& = & \left(\begin{array}{r}1\end{array}\right); & \vu_1\vone = 1
\end{array}
\]

Note that we have \(\vu_t \vone = {5 \choose t-1}\), as predicted by Remark~\ref{remark-sumcoeffs}.  Note also that the normalized version of \(\vu_3\) is
\[
\recip{\vu_3 \vone}\vu_3 = \recip{10} \left(\begin{array}{rrr}1 & 6 & 3\end{array}\right) =  \left(\begin{array}{rrr}\frac{1}{10} & \frac{3}{5} & \frac{3}{10}\end{array}\right) = L_3,
\]
and the same is true for the other \(\vu_t\)'s and \(L_t\)'s.  This elegant pattern for the landing vectors is the answer to Question~\ref{question-landvector}, but we need several more theorems to justify it.  The first two state that the \(\vu_t\) are eigenvectors for the probability transition matrices \(A_t\), and they are also ``chain eigenvectors'' in the sense that \(\vu_t A_{t,t-1}\) is a scalar multiple of \(\vu_{t-1}\):

\begin{thm} \label{theorem-samelevel}  
\[
\vu_t A_t = \vu_t d_t\mbox{, where } d_t := \frac{(n-t)(n+t-1)}{n(n-1)}.
\]
\end{thm}

We will discuss the proof of Theorem~\ref{theorem-samelevel} below.

\begin{cor} \label{corollary-eigen}\(\vu_t\) is a left eigenvector for the matrix \((I - A_t)\inv\) with eigenvalue \recip{1-d_t}.  \qed
\end{cor}

\begin{thm} \label{theorem-descend}
\[
\vu_t A_{t,t-1} = \vu_{t-1} h_t\mbox{, where } h_t := \frac{t(n-t+1)}{n(n-1)}.
\]
\end{thm}

Surprisingly, in our work later, we will only use the fact that \(\vu_t A_{t,t-1}\) is a multiple of \(\vu_{t-1}\); the actual value of \(h_t\) is immaterial.  We will explain this after Theorem~\ref{theorem-totaltime} below.

The proof of Theorem~\ref{theorem-samelevel} (respectively, Theorem~\ref{theorem-descend}) depends on some careful combinatorial bookkeeping.  We will suppress the computational details of the proofs, partly because of the tedium involved and partly because we have an independent way to answer Question~\ref{question-totaltime} that we will present in full detail later.  Instead, we will just give a sketch here and then illustrate with a numerical example.   

The main idea of both proofs is to track which states \(\vs \in S(n,t)\) (respectively, \(\vs \in S(n,t-1)\)) can be reached directly from which states \(\vr \in S(n,t)\), which we denote by \(\vr \rightarrow \vs\).  For \(\vr \rightarrow \vs\), we define \(\delta(\vr,\vs)\) to be the number of possible winner-loser pairs in state \vr\ that will take us to state \vs, that is, the numerator of the corresponding entry in \(A_t\) (respectively \(A_{t,t-1}\)), where the denominator is \(n(n-1)\).

The key step in the proof of Theorem~\ref{theorem-samelevel} is then to switch from summing over \vr\ to summing over \vs:
\qar
\vu_t A_t & = & \sum_{\vr \in S(n,t)} m_\vr \vr A_t\comment{by definition of \(\vu_t\)}\\
& = & \recip{n(n-1)} \sum_{\vr \in S(n,t)} \sum_{\{\vs: \vr \rightarrow \vs\}} m_\vr \delta(\vr,\vs) \vs\comment{by the action of \(A_t\)}\\
& = & \recip{n(n-1)} \sum_{\vs \in S(n,t)} \sum_{\{\vr: \vr \rightarrow \vs\}} m_\vr  \delta(\vr,\vs) \vs\comment{switching the summation}\\
& = & \recip{n(n-1)} \sum_{\vs \in S(n,t)} (n-t)(n+t-1)m_\vs \vs\comment{(see below)}\\
& = & \vu_t \frac{(n-t)(n+t-1)}{n(n-1)}\comment{by definition of \(\vu_t\)}
\raq
The work is in justifying the second to last equality above that 
\[
\sum_{\{\vr: \vr \rightarrow \vs\}} m_\vr  \delta(\vr,\vs) = t(n-t+1)m_\vs.
\]
This requires several pages of unenlightening calculation.  The proof of Theorem~\ref{theorem-descend} is similar, and similarly tedious.  We have spared you the full details, and instead we will illustrate with a larger concrete example.  Let \(n=10\) and \(t=4\); then in partition notation we have 
\qar
S(10,4) & = & \{3331,3322,4321,4411,4222,5311,5221,6211,7111\} \\
S(10,3) & = & \{433,442,541,532,631,622,721,811\}
\raq
and in vector notation we have
\qar
S(10,4) & = & \{(103),(022),(1111),(2002),(0301),(20101),(12001),(210001),(3000001)\}\\
S(10,3) & = & \{(0021),(0102),(10011),(01101),(101001),(020001),(1100001),(20000001)\}. 
\raq
Then 
\[
\vu_4 = \left(\begin{array}{rrrrrrrrr}4 & 6 & 24 & 6 & 4 & 12 & 12 & 12 & 4\end{array}\right), \vu_3 = \left(\begin{array}{rrrrrrrrr}3 & 3 & 6 & 6 & 6 & 3 & 6 & 3\end{array}\right)
\]
with \(\vu_4\vone = 84 = {9 \choose 3}\) and \(\vu_3\vone = 36 = {9 \choose 2}\), as predicted by Remark~\ref{remark-sumcoeffs}. 

The corresponding blocks of the transition matrix are 
\[
\left(\begin{array}{l|l} A_4 & A_{43}\end{array}\right) = \recip{90}\left(\begin{array}{rrrrrrrrr|rrrrrrrrr}
18 & 9 & 54 & 0 & 0 & 0 & 0 & 0 & 0 & & 9 & 0 & 0 & 0 & 0 & 0 & 0 & 0\\
8 & 40 & 24 & 0 & 18 & 0 & 0 & 0 & 0 & & 0 & 0 & 0 & 0 & 0 & 0 & 0 & 0\\
8 & 4 & 40 & 6 & 3 & 8 & 12 & 0 & 0 & & 2 & 3 & 0 & 4 & 0 & 0 & 0 & 0\\
0 & 0 & 16 & 24 & 0 & 32 & 0 & 0 & 0 & & 0 & 2 & 16 & 0 & 0 & 0 & 0 & 0\\
0 & 24 & 24 & 0 & 18 & 0 & 24 & 0 & 0 & & 0 & 0 & 0 & 0 & 0 & 0 & 0 & 0\\
0 & 0 & 10 & 15 & 0 & 26 & 6 & 15 & 0 & & 0 & 0 & 6 & 2 & 10 & 0 & 0 & 0\\
0 & 0 & 20 & 0 & 5 & 8 & 28 & 20 & 0 & & 0 & 0 & 0 & 4 & 0 & 5 & 0 & 0\\
0 & 0 & 0 & 0 & 0 & 12 & 12 & 36 & 12 & & 0 & 0 & 0 & 0 & 4 & 2 & 12 & 0\\
0 & 0 & 0 & 0 & 0 & 0 & 0 & 21 & 42 & & 0 & 0 & 0 & 0 & 0 & 0 & 6 & 21
\end{array}\right).
\]
Note that \(\vu_4 A_4 = \vu_4 \frac{13}{15} = \vu_4 d_4\) and \(\vu_4 A_{43} = \vu_3 \frac{14}{45} = \vu_3 h_4\), as predicted by Theorem~\ref{theorem-samelevel} and Theorem~\ref{theorem-descend}.  

\longcomment{
\bigskip

\textit{Proof of Theorem~\ref{theorem-descend}}.  For \(\vr \in S(n,t), \vs \in S(n,t-1)\), we say that {\em \vs\ follows \vr}, denoted \(\vr \rightarrow \vs\), if  it is possible to go directly from state \vr\ to state \vs\ by absorbing a solo player into a team.  This is reflected by a positive value of the corresponding entry in the matrix \(A_{t,t-1}\).  Note that \(\vr \rightarrow \vs\) if and only if \(\vs = \vr + \vb_i\) for some \(i \geq 1\), where we define
\[
\vb_i := \left(\begin{array}{rrrr}-1 & 0 & 0 & \cdots\end{array}\right) + \left(\begin{array}{rrrrrrr} 0 & \cdots & 0 & \underbrace{-1}_i & \underbrace{1}_{i+1} & 0 & \cdots\end{array}\right).
\]

For example, with \(n=10\) and \(t=4\) above, we have \((1111) \rightarrow (0021)\) (since \((0021) = (1111) + \vb_2\)) and \((1111) \rightarrow (0102)\) (since \((0102) = (1111) + \vb_3\)), but \((1111) \nrightarrow (10011)\).  

Now let \(\vr = \rowr\).  For \(\vs = \vr + \vb_i\), we define \(\delta(\vr,\vs)\) to be the number of possible winner-loser pairs in state \vr\ that will take us to state \vs, that is, the numerator of the corresponding entry in \(A_{t,t-1}\), where the denominator is \(n(n-1)\).  We now consider \vr\ and \vs\ as basis vectors in \(VS(n,t)\) and \(VS(n,t-1)\) respectively and \(A_{t,t-1}\) as a linear transformation from \(VS(n,t)\) to \(VS(n,t-1)\).  Then we have
\[
\vr A_{t,t-1} = \recip{n(n-1)}\sum_{\{\vs: \vr \rightarrow \vs\}} \delta(\vr,\vs) \vs.
\]
Note that if \(\vs = \vr + \vb_i\), then
\[
\delta(\vr,\vs) = \left\{\begin{array}{ll}r_1(r_1-1) & \mbox{if } i = 1, \\ r_1 r_i i & \mbox{if } i > 1.\end{array}\right.
\]
This is because there are \(r_1(r_1-1)\) ways to choose a winner and a loser that are both solo players, and there are \(r_1 r_i i\) ways to choose a winner on a team of size \(i > 1\) and a solo loser.

To prove Theorem~\ref{theorem-descend}, we switch from summing over \(\vr \in S(n,t)\) to summing over \(\vs \in S(n,t-1)\):
\qar
\vu_t A_{t,t-1} & = & \sum_{\vr \in S(n,t)} m_\vr \vr A_{t,t-1}\comment{by definition of \(\vu_t\)}\\
& = & \recip{n(n-1)} \sum_{\vr \in S(n,t)} \sum_{\{\vs: \vr \rightarrow \vs\}} m_\vr \delta(\vr,\vs) \vs\\
& = & \recip{n(n-1)} \sum_{\vs \in S(n,t-1)} \sum_{\{\vr: \vr \rightarrow \vs\}} m_\vr  \delta(\vr,\vs) \vs
\raq
It now suffices to show that for each \(\vs \in S(n,t-1)\), we have
\[
\sum_{\{\vr: \vr \rightarrow \vs\}} m_\vr  \delta(\vr,\vs) = t(n-t+1)m_\vs.
\]
Suppose \(\vs = \rows\).  The \vr's that are followed by \vs\ are those for which \(\vs = \vr + \vb_i\) for \(1 \leq i \leq k-1\).  We define \(\vr_i := \vs - \vb_i\).  Then if \(i = 1\), we have 
\qar
\vr_1 & = & \vs - \vb_1 = \left(\begin{array}{rrrrr}s_1+2 & s_2 - 1 & s_3 & \cdots & s_k\end{array}\right)\\
\delta(\vr_1,\vs) & = & (s_1+2)(s_1+1)\\
m_{\vr_1}\delta(\vr_1,\vs) & = & \frac{t!(s_1+2)(s_1+1)}{(s_1+2)!(s_2-1)!\cdots s_k!}\\
& = & \frac{ts_2(t-1)!}{s_1!s_2!\cdots s_k!}\\
& = & ts_2m_\vs.
\raq
If \(i > 1\), we have 
\qar
\vr_i & = & \vs - \vb_i = \left(\begin{array}{rrrrrrr}s_1+1 & s_2 & \cdots & s_i+1 & s_{i+1} - 1 & \cdots & s_k\end{array}\right)\\
\delta(\vr_i,\vs) & = & (s_1+1)(s_i+1)i\\
m_{\vr_i}\delta(\vr_i,\vs) & = & \frac{t!(s_1+1)(s_i+1)i}{(s_1+1)!s_2!\cdots (s_i+1)!(s_{i+1} - 1)! \cdots s_k!}\\
& = & \frac{ti s_{i+1}(t-1)!}{s_1!s_2!\cdots s_k!}\\
& = & ti s_{i+1} m_\vs.
\raq
We now calculate:
\qar
\sum_{\{\vr: \vr \rightarrow \vs\}} m_\vr  \delta(\vr,\vs) & = & \sum_{i=1}^{k-1} m_{\vr_i}\delta(\vr_i,\vs)\\
& = & m_{\vr_1}\delta(\vr_1,\vs) + \sum_{i=2}^{k-1} m_{\vr_i}\delta(\vr_i,\vs)\\
& = & ts_2m_\vs + \sum_{i=2}^{k-1} ti s_{i+1} m_\vs\\
& = & t\left[s_2 + 2s_3 + \cdots (k-1)s_k\right]m_\vs\\
& = & t\left[(s_1 + 2s_2 + 3s_3 + \cdots ks_k) - (s_1 + \cdots + s_k)\right]m_\vs\\
& = & t[n-(t-1)]m_\vs\comment{since \(\vs \in S(n,t-1)\)}\\
& = &t(n-t+1)m_\vs
\raq
This is what we needed to prove.  \qed
\bigskip

\textit{Proof of Theorem~\ref{theorem-samelevel}}.  The proof is similar to that of Theorem~\ref{theorem-descend}, but it is a little trickier to figure out which partitions in \(S(n,t)\) follow which.  Suppose \(\vr,\vs \in S(n,t)\) with \(\vr = \rowr\), \(\vs = \rows\).  As in the proof of Theorem~\ref{theorem-descend}, we say {\em \vs\ follows \vr}, denoted \(\vr \rightarrow \vs\), if it is possible to go directly from state \vr\ to state \vs, in other words, if the corresponding entry in the matrix \(A_t\) is positive.  

For example, with \(n=10\) and \(t=4\) above, we have \((022) \rightarrow (103), (022) \rightarrow (022)\), and \((022) \rightarrow (1111)\), but \((022) \nrightarrow (2002)\).  

For \(\vr \rightarrow \vs\), we define \(\delta(\vr,\vs)\) as above to be the number of possible winner-loser pairs in state \vr\ that will take us to state \vs, that is, the numerator of the corresponding entry in \(A_t\), where the denominator is \(n(n-1)\).

As with Theorem~\ref{theorem-descend}, we prove Theorem~\ref{theorem-samelevel} by switching from summing over \vr\ to summing over \vs:
\qar
\vu_t A_t & = & \sum_{\vr \in S(n,t)} m_\vr \vr A_t\comment{by definition of \(\vu_t\)}\\
& = & \recip{n(n-1)} \sum_{\vr \in S(n,t)} \sum_{\{\vs: \vr \rightarrow \vs\}} m_\vr \delta(\vr,\vs) \vs\\
& = & \recip{n(n-1)} \sum_{\vs \in S(n,t)} \sum_{\{\vr: \vr \rightarrow \vs\}} m_\vr  \delta(\vr,\vs) \vs
\raq
It now suffices to show that for each \(\vs \in S(n,t)\), we have
\[
\sum_{\{\vr: \vr \rightarrow \vs\}} m_\vr  \delta(\vr,\vs) = (n-t)(n+t-1)m_\vs.
\]
In fact, we will show that for each \(\vs \in S(n,t)\), we have
\[
\sum_{\{\vr: \vr \rightarrow \vs\}} \frac{m_\vr}{m_\vs}  \delta(\vr,\vs) = (n-t)(n+t-1).
\]
We will list the various \vr\ for which \(\vr \rightarrow \vs\) and for each one, find their contribution to 
\[
\sum_{\{\vr: \vr \rightarrow \vs\}} \frac{m_\vr}{m_\vs}  \delta(\vr,\vs).\hskip 5mm \mbox{(*)}
\]
To list the ways in which \vs\ can follow \vr, we introduce some notation.  If the loser is from a group of size \(i > 1\) and the winner is from a group of size \(j \geq 1\) (a different group in the case \(i=j\)), we will have \(\vs = \vr + \va_i + \vb_j\), where we define
\qar
\va_i & := & \left(\begin{array}{rrrrrrr}0 & \cdots & 0 & \underbrace{1}_{i-1} & \underbrace{-1}_i & 0 & \cdots\end{array}\right)\\
\vb_j & := & \left(\begin{array}{rrrrrrr}0 & \cdots & 0 & \underbrace{-1}_j & \underbrace{1}_{j+1} & 0 & \cdots\end{array}\right).
\raq
Since we are listing the \vr's for a fixed \vs\ below, we will write this as \mbox{\(\vr = \vs - \va_i - \vb_j\)}.  Note that the nonzero entries of \(\va_i\) and \(\vb_j\) will overlap if \(i-2 \leq j \leq i\) and will be independent if \(j < i-2\) or \(j > i\).

Now, there are four ways in which we can have \(\vr \rightarrow \vs\):

\casI
\item  \(\vr = \vs\).  This can occur in two possible ways:
\num
\item  The loser and the winner are already in the same group, say of size \(i > 1\).  The contribution of this case to \(\delta(\vr,\vs)\) is 
\[
\sum_{i=2}^k i (i-1) r_i = \sum_{i=2}^k i(i-1) s_i.
\]
\item  The loser comes from a group of size \(i>1\) and the winner from a group of size \(i-1\).  This is the case \(\vr = \vs - \va_i - \vb_{i-1}\).  The contribution of this case to \(\delta(\vr,\vs)\) is 
\[
\sum_{i=2}^k i(i-1) r_i r_{i-1} = \sum_{i=2}^k i(i-1) s_i s_{i-1}.  
\]
\mun
Adding these two, we get
\begin{equation}
\delta(\vr,\vs) = \sum_{i=2}^k i(i-1) s_i (s_{i-1}+1).\label{net-one}
\end{equation}
In this case, \(\frac{m_\vr}{m_\vs} = 1\), so the contribution to (*) is just (1) above.  

\item  \(\vr = \vs - \va_i - \vb_i\) for some \(i > 1\).  This means that the loser and winner come from two different groups of size \(i\).  Then \(r_1 = s_1, \dots, r_{i-1} = s_{i-1} - 1, r_i = s_i + 2, r_{i+1} = s_{i+1}-1,\dots, r_k = s_k\).  Then
\qar
\delta(\vr,\vs) & = & i^2 r_i (r_i-1) = i^2 (s_i+2)(s_i+1)\\
\frac{m_\vr}{m_\vs} & = & \frac{s_1!\cdots s_k!}{s_1!\cdots (s_{i-1}-1)!(s_i+2)!(s_{i+1}-1)! \cdots s_k!}\\
\frac{m_\vr}{m_\vs}\delta(\vr,\vs) & = & i^2 s_{i-1}s_{i+1}
\raq
Hence, the contribution to (*) from this case is
\begin{equation}
\sum_{i=2}^{k-1} i^2 s_{i-1}s_{i+1}.\label{net-two}
\end{equation}

\item  \(\vr = \vs - \va_i - \vb_{i-2}\) for some \(i > 2\).  This means that the loser comes from a group of size \(i\) (which could be as large as \(k+1\), since after the loss the group would only have size \(k\)) and the winner from a group of size \(i-2\).    Then \(r_1 = s_1, \dots, r_{i-2} = s_{i-2} + 1, r_{i-1} = s_{i-1} - 2, r_i = s_i+1,\dots, r_k = s_k\).  Then
\qar
\delta(\vr,\vs) & = & i(i-2) r_i r_{i-2} = i(i-2) (s_i+1)(s_{i-2}+1)\\
\frac{m_\vr}{m_\vs} & = & \frac{s_1!\cdots s_k!}{s_1!\cdots (s_{i-2}+1)!(s_{i-1}-2)!(s_i+1)! \cdots s_k!}\\
\frac{m_\vr}{m_\vs}\delta(\vr,\vs) & = & i(i-2) s_{i-1}(s_{i-1}-1)
\raq
Hence, the contribution to (*) from this case is
\[
\sum_{i=3}^{k+1} i(i-2) s_{i-1}(s_{i-1}-1),
\]
which we re-index as
\begin{equation}
\sum_{i=2}^k (i+1)(i-1) s_i (s_i-1).\label{net-three}
\end{equation}

\item  \(\vr = \vs - \va_i - \vb_j\) for some \(i > 1\) and \(j < i-2\) or \(j > i\).  The loser comes from a group of size \(i\) (again, possibly as large as \(k+1\)) and the winner from a group of size \(j\).    Then \(r_1 = s_1, \dots, r_{i-1} = s_{i-1} - 1, r_i = s_i + 1, \dots, r_j = s_j+1,r_{j+1} = s_{j+1}-1,\dots, r_k = s_k\).  Then
\qar
\delta(\vr,\vs) & = & ij r_i r_j = ij (s_i+1)(s_j+1)\\
\frac{m_\vr}{m_\vs} & = & \frac{s_1!\cdots s_k!}{s_1!\cdots (s_{i-1}-1)!(s_i+1)!\cdots (s_j+1)!(s_{j+1}-1)! \cdots s_k!}\\
\frac{m_\vr}{m_\vs}\delta(\vr,\vs) & = & ij s_{i-1}s_{j+1}
\raq
Hence, the contribution to (*) from this case is
\begin{equation}
\underbrace{\sum_{i = 2}^{k+1} \sum_{j=1}^{k-1}}_{\{i,j:j<i-2\mbox{ or }j>i\}} ij s_{i-1}s_{j+1}.\label{net-four}
\end{equation}
\sacI
To finish the proof of Theorem~\ref{theorem-samelevel}, then, we must show that the sum of the contributions from (1), (2), (3), and (4), is \((n-t)(n+t-1)\).  We first expand:
\qar
(n-t)(n+t-1) & = & \left(\sum_{i=1}^k i s_i - \sum_{i=1}^k s_i\right)\left(\sum_{j=1}^k j s_j + \sum_{j=1}^k s_j - 1\right)\\
& = & \left(\sum_{i=1}^k (i-1) s_i\right)\left(\sum_{j=1}^k (j+1) s_j - 1\right)
\raq
The linear terms here match the linear terms from (1) and (3) above:
\[
\underbrace{-\sum_{i=1}^k (i-1)s_i}_{\mbox{from }(n-t)(n+t-1)} = \underbrace{\sum_{i=2}^k i(i-1)s_i}_{\mbox{from (1)}} - \underbrace{\sum_{i=2}^k (i+1)(i-1)s_i}_{\mbox{from (3)}}
\]
So now we need only confirm that the quadratic terms
\begin{equation}
\left(\sum_{i=1}^k (i-1) s_i\right)\left(\sum_{j=1}^k (j+1) s_j\right)
\end{equation}
equal the remaining contributions from (1), (2), (3), and (4):
\[
\sum_{i=2}^k i(i-1) s_i s_{i-1} + \sum_{i=2}^{k-1} i^2 s_{i-1}s_{i+1} + \sum_{i=2}^k (i+1)(i-1) s_i^2 + \underbrace{\sum_{i = 2}^{k+1} \sum_{j=1}^{k-1}}_{\{i,j:j<i-2\mbox{ or }j>i\}} ij s_{i-1}s_{j+1}
\]
But this is now clear, since the terms in (5) with \(i=j\) are equal to the third sum above, those with \(i = j+1\) give the first, those with \(i = j+2\) give the second, and those with \(i<j\) and \(i > j+2\) give the fourth.  This completes the proof of Theorem~\ref{theorem-samelevel}.  \qed
}

We can now justify our answer to Question~\ref{question-landvector}:

\begin{thm} \label{theorem-landing}For all \(t\), the landing vector \(L_t\) is the normalized \(\vu_t\), that is,
\[
L_t = \recip{\vu_t \vone} \vu_t.
\]
\end{thm}

\textit{Proof}.  First, we note that \(L_n = (1) = \vu_n\).  Proceeding downwards by induction, we assume the theorem for \(L_t\) and show it for \(L_{t-1}\):  
\qar
L_{t-1} & = & L_tP_{t,t-1}\comment{by construction of \(L_{t-1}\)}\\
& = & L_t (I - A_t)\inv A_{t,t-1}\comment{by construction of \(P_{t,t-1}\)}\\
& = & \recip{\vu_t\vone}\vu_t(I - A_t)\inv A_{t,t-1}\comment{by the induction hypothesis}\\
& = & \recip{\vu_t \vone} \recip{1-d_t} \vu_t A_{t,t-1}\comment{by Corollary~\ref{corollary-eigen}}\\
& = & \recip{\vu_t \vone} \recip{1-d_t} \vu_{t-1} h_t\comment{by Theorem~\ref{theorem-descend}}
\raq
This shows that \(L_{t-1}\) is a scalar multiple of \(\vu_{t-1}\).  But since we know that \(L_{t-1}\) is a probability vector, i.e., that its entries sum to one, we must have that 
\[
L_{t-1} = \recip{\vu_{t-1} \vone} \vu_{t-1},
\]
as desired.  \qed

\begin{rem} \label{remark-findcoeffs}The proof of Theorem~\ref{theorem-landing} gives an alternate way to find \(\vu_t \vone\). \end{rem}
\textit{Proof}.  We can find a relationship between \(\vu_{t-1}\vone\) and \(\vu_t\vone\):
\qar
L_{t-1} \vone & = & \frac{h_t}{1-d_t} \recip{\vu_t \vone} \vu_{t-1} \vone\comment{from the proof above}\\
1 & = & \frac{h_t}{1-d_t} \recip{\vu_t \vone} \vu_{t-1} \vone\comment{since \(L_{t-1}\) is a probability vector}\\
\vu_{t-1}\vone & = & \frac{1-d_t}{h_t}\vu_t\vone\comment{by cross multiplication}\\
& = & \frac{t-1}{n-t+1} \vu_t\vone\comment{by definition of \(d_t\) and \(h_t\)}
\raq
This gives us the recursive system
\qar
\vu_n \vone & = & 1\\
\vu_{n-1} \vone & = & \frac{n-1}1 \vu_n \vone = \frac{n-1}1\\
\vu_{n-2} \vone & = & \frac{n-2}2 \vu_{n-1} \vone = \frac{n-2}2 \frac{n-1}1\\
& \vdots & \\
\vu_t \vone & = & \frac t{n-t} \cdots \frac{n-2}2 \frac{n-1}1 = {n-1 \choose t-1},
\raq
confirming our result from Remark~\ref{remark-sumcoeffs}. \qed

\section{Expected times}

We are now ready to answer Questions~\ref{question-stagetime} and~\ref{question-totaltime}.  Recall that \(L_t\) is the row vector whose \(j\)-th entry is the probability that our first arrival in Stage \(t-1\) from Stage \(t\) is in the \(j\)-th state in Stage \(t-1\).  We define \(e_{t,t-1}\) to be the expected time from our first arrival in Stage \(t\) to our first arrival in Stage \(t-1\).  We have an immediate answer for Question~\ref{question-stagetime}.

\begin{thm} \label{theorem-stagetime} 
\[
e_{t,t-1} = \frac{n(n-1)}{t(t-1)} = \frac{{n \choose 2}}{{t \choose 2}}
\]
\end{thm}

\textit{Proof}.  When we worked out the case for \(n=4\) we derived a formula 
for \(e_{t,t-1}\) that clearly generalizes to larger cases.  (This is a standard result in the theory of Markov chains; see Theorem 3.3.5 in \cite{kemeny}.)  We proceed from that formula:
\qar
e_{t,t-1} & = & L_t(I-A_t)\inv \vone\\
& = & L_t \recip{1-d_t} \vone\comment{by Corollary~\ref{corollary-eigen} and Theorem~\ref{theorem-landing}}\\
& = & \recip{1-d_t} L_t \vone\comment{since \recip{1-d_t} is a scalar}\\
& = & \recip{1-d_t}\comment{since \(L_t\) is a probability vector}\\
& = & \frac{n(n-1)}{t(t-1)}\comment{by definition of \(d_t\)}
\raq
\qed

We now just add the times at each stage to answer Question~\ref{question-totaltime}:

\begin{thm} \label{theorem-totaltime}The expected time to the final state is \((n-1)^2\).
\end{thm}

\textit{Proof}.  We use a partial fraction expansion:
\qar
\sum_{t=2}^n e_{t,t-1} & = & \sum_{t=2}^n \frac{n(n-1)}{t(t-1)}\comment{by Theorem~\ref{theorem-stagetime}}\\
& = & n(n-1) \sum_{t=2}^n \point{\recip{t-1}-\recip{t}},\comment{a telescoping series}\\
& = & n(n-1)\point{1-\recip{n}}\\
& = & (n-1)^2
\raq
\qed

One slightly surprising element of the proofs above is that we never used the formula for the ``chain eigenvalue'' \(h_t\) from Theorem~\ref{theorem-descend}.  (We did use the value of \(h_t\) in the proof of Remark~\ref{remark-findcoeffs}, but Remark~\ref{remark-findcoeffs} was not used to prove anything else.)  This is less surprising when we realize that the value of \(h_t\) can be derived from the value of \(d_t\) by the following method, which is independent of the formula in Theorem~\ref{theorem-descend}.  Note that the row vector \(L_t\left(\begin{array}{ll} A_t & A_{t,t-1}\end{array}\right)\) gives the complete set of probabilities of landing in the various states in Stage \(t\) and Stage \(t-1\) one step after landing in Stage \(t\).  As such, the entries this row vector add to one.  But we can calculate this vector:
\qar
L_t A_t & = & L_t d_t \comment{by Theorems~\ref{theorem-samelevel} and~\ref{theorem-landing}}\\
L_t A_{t,t-1} & = & \vu_t \recip{\vu_t \vone}A_{t,t-1}\comment{by Theorem~\ref{theorem-landing}}\\
& = &  \vu_{t-1}\frac{h_t}{\vu_t \vone}\comment{by Theorem~\ref{theorem-descend}}\\
& = &  L_{t-1}\frac{(\vu_{t-1} \vone) h_t}{\vu_t \vone}\comment{by Theorem~\ref{theorem-landing}}\\
\raq
Therefore,
\qar
L_t\left(\begin{array}{ll} A_t & A_{t,t-1}\end{array}\right)\vone & = & 1\comment{by the discussion above}\\
\left(\begin{array}{ll} L_t d_t & L_{t-1}\frac{(\vu_{t-1} \vone) h_t}{\vu_t \vone}\end{array}\right)\vone & = & 1\comment{by the calculations immediately above}\\
d_t + \frac{(\vu_{t-1} \vone) h_t}{\vu_t \vone} & = & 1\comment{since \(L_t\) and \(L_{t-1}\) are probability vectors}\\
d_t + \frac{(t-1) h_t}{n-t+1} & = & 1\comment{by Remark~\ref{remark-sumcoeffs}.}\\
\raq
Thus, \(d_t\) and \(h_t\) are dependent on each other, and if we use a particular value of one, then we are also implicitly using the corresponding value of the other.  And note that the value of \(d_t\) did indeed play a key role in the proof of Theorem~\ref{theorem-stagetime} above. 

\section{A symmetric approach}

Although we think the answers to Questions~\ref{question-stagetime} and~\ref{question-landvector} are interesting in their own right, we can derive the answer to Question~\ref{question-totaltime} independently without going through the calculations above.  In particular, this method does not rely on the omitted proofs of Theorems~\ref{theorem-samelevel} and~\ref{theorem-descend}.

We start with \(n\) players, each of whom initially represents a different field. We arbitrarily choose one field to focus on, say, statistics.  At any point in the game, we define a set of random variables \(x_0, \dots, x_n\), where \(x_i\) represents the number of {\em {future wins by statisticians}}, given that there are \(i\) statisticians currently remaining.  (Note that it does not matter what the configuration of the other \(n-i\) players into teams is.)  We have easy boundary values:  \(x_0 = 0\), since if statistics has been wiped out as a field, then there can be no future converts to statistics; and \(x_n = 0\), since if everyone is now a statistician then the game is over.  

We now set up a system of equations for the other \(x_i, 1 \leq i \leq n-1\).  In each round, there are \(n(n-1)\) choices for the winner and loser.  With \(i\) statisticians currently, there are four possibilities for how the winner and loser can be arranged with respect to the statisticians:
\num
\item  Both winner and loser are statisticians.  There are \(i(i-1)\) ways this can happen.  The number of wins by statisticians has gone up by one, and the new expectation at the following round is again \(x_i\), since we again have \(i\) statisticians.  

\item  Only the winner is a statistician.  There are \(i(n-i)\) ways this can happen.  The number of wins by statisticians has gone up by one, and the new expectation at the following round is \(x_{i+1}\), since we then have \(i+1\) statisticians.  

\item  Only the loser is a statistician.  There are \(i(n-i)\) ways this can happen.  The number of wins by statisticians is unchanged, and the new expectation at the following round is \(x_{i-1}\)
.  

\item  Neither the winner nor the loser is a statistician.  There are \mbox{\((n-i)(n-i-1)\)} ways this can happen.  The number of wins by statisticians is unchanged, and the new expectation at the following round is again \(x_i\)
.  
\mun
This gives us the following equation:
\[
x_i = \frac{i(i-1)}{n(n-1)}(1+x_i) +\frac{ i(n-i)}{n(n-1)}(1+x_{i+1}) + \frac{i(n-i)}{n(n-1)}x_{i-1} + \frac{(n-i)(n-i-1)}{n(n-1)}x_i
\]
Mercifully, this simplifies rather dramatically:
\[
2x_i - (x_{i-1} + x_{i+1}) = \frac{n-1}{n-i}
\]
This gives us a linear system for the \(x_i\)'s:
\[
\left(\begin{array}{rrrrr}
2 & -1 & 0 & \cdots & 0\\
-1 & 2 & -1 & \cdots & 0\\
0 & -1 & 2 & \cdots & 0\\
\vdots & \vdots & \vdots & \ddots & \vdots\\
0 & 0 & 0 & \cdots & 2\end{array}\right)
\left(\begin{array}{c}
x_1\\
x_2\\
x_3\\
\vdots\\
x_{n-1}
\end{array}\right) = 
\longcomment{
\left(\begin{array}{rrrrr}
\recip{n-1} &  & &  & \\
 & \recip{n-2} &  & & \\
&  & \recip{n-3} &  & \\
& & & \ddots & \\
& & &  & 1\end{array}\right)
\left(\begin{array}{c}
n-1\\
n-1\\
n-1\\
\vdots\\
n-1\\
\end{array}\right)
}
\left(\begin{array}{c}
1\\
\frac{n-1}{n-2}\\
\frac{n-1}{n-3}\\
\vdots\\
n-1\\
\end{array}\right)
\]
We denote the \((n-1)\times(n-1)\) matrix on the left by \(M_n\).  It is an amusing exercise to compute \(M_n\inv\); for example, with \(n=6\) we have
\[
M_6 = \left(\begin{array}{rrrrr}
2 & -1 & &  & \\
-1 & 2 & -1 & & \\
&  -1 & 2 & -1 & \\
& & -1 & 2 & -1\\
& & &  -1 & 2\end{array}\right), 
M_6\inv = \recip{6}\left(\begin{array}{rrrrr}
5 & 4 & 3 & 2 & 1\\
4 & 8 & 6 & 4 & 2\\
3 & 6 & 9 & 6 & 3\\
2 & 4 & 6 & 8 & 4\\
1 & 2 & 3 & 4 & 5\end{array}\right).
\]
The pattern in the right-hand matrix is that the \((i,j)\)-entry is \(i(n-j)\) for entries above the main diagonal and \(j(n-i)\) for entries below.  In other words,
\[
\point{M_n\inv}_{i,j} = \recip{n}\min\{i,j\}\left[n-\max\{i,j\}\right].
\]
To answer Question~\ref{question-totaltime}, we need to know the expected number of future wins by statisticians at the very start of the game.  We start with one statistician, so we solve our system for \(x_1\) using the first row of \(M_n\inv\):
\qar
\left(\begin{array}{c}
x_1\\
x_2\\
\vdots\\
x_{n-1}
\end{array}\right) & = &
M_n\inv
\left(\begin{array}{c}
1\\
\frac{n-1}{n-2}\\
\vdots\\
n-1\\
\end{array}\right)\\
x_1 & = & \recip{n}\left(\begin{array}{cccc}n-1 & n-2 & \cdots & 1\end{array}\right)
\left(\begin{array}{cccc}1 & \frac{n-1}{n-2} & \cdots & n-1\end{array}\right)^T\\
& = & \recip n\left[(n-1) + (n-1) +\cdots +(n-1)\right]\\
& = & \frac{(n-1)^2}n
\raq
We have just computed the expected number of total wins by {\em {statisticians}}.  By symmetry, every other field expects the same number of wins, so the total number of rounds of the game (again, exploiting linearity of expectation) is \(n\frac{(n-1)^2}n = (n-1)^2\).  This confirms our answer to Question~\ref{question-totaltime} from the small games and the derivation in the previous section.  

Finally, we address the temptation to hope that a Markov chain with such a nice expectation might also have an interesting variance.  Following Theorem 3.3.5 in \cite{kemeny}, we can compute the variance of the time to absorbtion via the matrix \mbox{\(N := (I-A)\inv\)}, where \(A\) is the submatrix of \(P\) obtained by deleting the final row and column, which correspond to the absorbing state.  We then define the column vector \(\tau := N\vone\) (the expected time to absorbtion from each state), and let \(\tau_{\mbox{\tiny sq}}\) be the column vector whose entries are the squares of those in \(\tau\).  Then \cite{kemeny} tells us that the variance of the time to absorbtion is the first entry of the vector
\[
\tau_2 := (2N-I)\tau - \tau_{\mbox{\tiny sq}}.
\]
For \(n = 2,3,4\), the variances turn out to be 0,6, and 32, raising the hope that an interesting sequence of integers might ensue.  Sadly, for \(n=5\) and \(n=6\), the variances are \(\frac{890}9\) and \(\frac{469}2\), respectively.  We challenge you to discover, prove, and interpret the general pattern!\\

\small
\noindent {\bf Acknowledgement}  We thank John Brevik for suggesting this problem and Kent Merryfield and Peter Ralph for useful conversations.

\small

\noindent {\bf Summary}  Consider a system of \(n\) players in which each initially starts on a different team.  At each time step, we select an individual winner and an individual loser randomly and the loser joins the winner's team.  The resulting Markov chain and stochastic matrix clearly have one absorbing state, in which all players are on the same team, but the combinatorics along the way are surprisingly elegant.  The expected number of time steps until each team is eliminated is a ratio of binomial coefficients.  When a team is eliminated, the probabilities that the players are configured in various partitions of \(n\) into \(t\) teams are given by multinomial coefficients.  The expected value of the time to absorbtion is \((n-1)^2\) steps.  The results depend on elementary combinatorics, linear algebra, and the theory of Markov chains.  \\

\noindent {\bf  ROBERT MENA} joined the faculty at Long Beach State in 1988 after 15 years in the faculty at the University of Wyoming, and graduate school at the University of Houston. He has also spent time at Caltech and Ohio State. He is an enthusiastic solver of quote and other sorts of acrostics as well as ken ken puzzles. Even after 40 years of teaching, he still purports to enjoy the company of his students and the thrill of teaching mathematics.\\

\noindent {\bf WILL MURRAY} did graduate work in algebra and is also interested in probability, combinatorics, and analysis.  He has taught at Berkeley, Long Beach State, the Royal University of Phnom Penh, and the University of M\'{e}d\'{e}a in Algeria.  An avid juggler and traveler, he has performed on five continents and written and lectured about the mathematics of juggling.  At home, he enjoys playing with his pets, one of whom starred on {\em {America's Got Talent}}.

\vfill\eject
\end{document}